\documentclass[conference,a4paper,10pt]{IEEEtran}
 \pdfoutput=1
\usepackage{amsmath}
\usepackage{amsthm}
\usepackage{amssymb}
\usepackage{stmaryrd}
\usepackage[usenames,dvipsnames]{xcolor}

\usepackage[noadjust]{cite}

\usepackage{hyperref}
\hypersetup{
    bookmarks=true,         
    unicode=false,          
    pdftoolbar=true,        
    pdfmenubar=true,        
    pdffitwindow=false,     
    pdfstartview={FitH},    
    pdfnewwindow=true,      
    colorlinks=true,       
    linkcolor=blue!80!black,          
    citecolor=Maroon,        
    filecolor=Maroon,      
    urlcolor=OliveGreen           
}

\usepackage{multirow}
\usepackage{caption2}
\usepackage{longtable}
\usepackage[all,color, pdftex]{xy}
\usepackage{enumerate}
\usepackage{cases}

\newtheorem{theorem}{Theorem}

\newtheorem{lemma}[theorem]{Lemma}
\newtheorem{proposition}[theorem]{Proposition}

\theoremstyle{remark}
\newtheorem{example}{Example}

\newcommand{\bz}{{\mathbb Z}}
\newcommand{\su}{{\sf u}}
\newcommand{\sv}{{\sf v}}

\newcommand{\C}{\mathcal C}

\newcommand{\F}{\mathcal F}
\newcommand{\G}{\mathcal G}
\newcommand{\HH}{\mathcal H}

\newcommand{\w}{\overline{w}}
\newcommand{\ind}{\textrm{in}}
\newcommand{\outd}{\textrm{out}}
\newcommand{\supp}{{\rm supp}}
\newcommand{\etal}[1]{{\it et al.}#1}
\newcommand{\floor}[1]{\left\lfloor#1\right\rfloor}

\begin{document}

\sloppy

\title{Decompositions of Edge-Colored Digraphs: A New Technique in the Construction of Constant-Weight Codes and Related Families}

\author{\IEEEauthorblockN{Yeow Meng Chee\IEEEauthorrefmark{1}, 
Fei Gao\IEEEauthorrefmark{2},
Han Mao Kiah\IEEEauthorrefmark{1}, 
Alan Chi Hung Ling\IEEEauthorrefmark{3}, 
Hui Zhang\IEEEauthorrefmark{1}, 
Xiande Zhang\IEEEauthorrefmark{1}} 
\IEEEauthorblockA{\IEEEauthorrefmark{1}
\small School of Physical and Mathematical Sciences, Nanyang Technological University, Singapore,\\
emails: ymchee@ntu.edu.sg, hmkiah@ntu.edu.sg, huizhang@ntu.edu.sg, xiandezhang@ntu.edu.sg} 
\IEEEauthorblockA{\IEEEauthorrefmark{2}\small Institute of High Performance Computing, Agency for Science, Technology and Research, Singapore,
email: gaofei@ihpc.a-star.edu.sg} 
\IEEEauthorblockA{\IEEEauthorrefmark{3}\small Department of Computer Science, University of Vermont, USA, 
email: aling@emba.uvm.edu}
 \vspace*{-0.2in}
}

\date{}
\maketitle

\begin{abstract}

We demonstrate that certain Johnson-type bounds are asymptotically exact
for a variety of classes of codes, namely, constant-composition codes,
nonbinary constant-weight codes and multiply constant-weight codes.
This was achieved via an interesting application of the theory of 
decomposition of edge-colored digraphs.

 \vspace*{-0.1in}
%
%
\end{abstract}

\section{Introduction}

In 1970s, Wilson \cite{Wilson:1972a,Wilson:1972b,WIlson:1972d,Wilson:1975} 
demonstrated that the elementary necessary conditions for the 
existence of balanced incomplete block designs are asymptotically sufficient.
His work therefore determined the size of an optimal binary constant-weight code
of weight $w$ and distance $2w-2$, 
provided that the length $n$ is sufficiently large and satisfies certain congruence classes.
Moreover, these optimal codes meet the Johnson bound
and his results show that Johnson bound is asymptotically exact for fixed weight $w$ and distance $2w-2$.

Since Wilson's seminal work, there have been interesting developments
in both areas of combinatorial design and coding theory. 
In the former, Wilson's ideas matured into Lamken and Wilson's theory of {decomposition of edge-colored digraphs} \cite{LamkenWilson:2000}
and the theory has been used extensively in establishing the asymptotic existence 
of many classes of combinatorial designs. 
Theoretical developments also extend Lamken and Wilson's results to other classes of decompositions.
Of particular interest is the class of {superpure} decompositions \cite{Hartmann:2002}.

On the other hand, binary constant-weight codes have been generalized to 
constant-composition codes and nonbinary constant-weight codes 
and have found applications in
powerline communications \cite{Chuetal:2004,Cheeetal:2013},
frequency hopping \cite{Chuetal:2006},
coding for bandwidth-limited channels \cite{CostelloForney:2007}. 
More recently, 
multiply constant-weight codes are introduced 
in an application for physically unclonable functions \cite{Cheeetal:preprint2013}.
Not surprisingly, in these generalizations, Johnson-type upper bounds 
have been derived.

Therefore, a natural question is whether the advanced techniques of decompositions of edge-colored digraphs
are relevant in constructing optimal codes  
and whether Johnson-type bounds are asymptotically exact for these generalizations.
We answer both questions in the affirmative for certain distances 
where the weight is fixed.

\section{Preliminary}


The ring $\bz/q\bz$ is denoted by $\bz_q$. 
For positive integer $n$, the set $\{1,2,\ldots,n\}$ is denoted by $[n]$.


Let $\bz_q^X$ denoted the set of vectors whose elements belong to $\bz_q$
and are indexed by $X$.
A {\em $q$-ary code of length $n$} is then a set $\C\subseteq\bz_q^X$ with $|X|=n$.
 The {\em support} of a vector $\su\in \bz_q^X$, denoted $\supp(\su)$, is the set $\{x\in X:\su_x\neq 0\}$. 
The {\em Hamming weight} of $\su\in \bz_q^X$ is defined as $||\su||=|\supp(\su)|$. 
The distance induced by this norm is called the {\em Hamming distance}, so that $d(\su,\sv)=||\su-\sv||$, for $\su,\sv\in \bz_q^X$. A code $\C$ is said to have {\em distance $d$} if $d(\su,\sv)\geq d$ for all distinct $\su,\sv\in \C$. 
The {\em composition} of a vector $\su\in \bz_q^X$ is the tuple $\w=[w_1,\dots,w_{q-1}]$, 
where $w_i=|\{x\in X:\su_x=i\}|$, where $i\in \bz_q\setminus\{0\}$. 
Unless mentioned otherwise, we always assume $w_1\ge w_2\ge \cdots\ge w_{q-1}$.

\subsection{Constant-Weight Codes and Constant-Composition Codes}
 A code $\C$ is said to have {\em constant-weight $w$} if every codeword in $\C$ has weight $w$, and 
 has {\em constant-composition $\w$} if every codeword in $\C$ has composition $\w$. 
 We refer to a $q$-ary code of length $n$, distance $d$, and constant-weight $w$ as a CWC$(n,d,w)_q$. 
 If in addition, the code has constant-composition $\w$, then it is referred to as a CCC$(n,d,\w)_q$.
  The maximum size of a CWC$(n,d,w)_q$ is denoted $A_q(n,d,w)$ while 
  the maximum size of a CCC$(n,d,\w)_q$ is denoted $A_q(n,d,\w)$. 
  Any CWC$(n,d,w)_q$ or CCC$(n,d,\w)_q$ attaining the maximum size is called {\em optimal}.

Consider a composition $\w$ and let $w=\sum_{i=1}^{q-1}w_i$.
Then a CCC$(n,d,\w)_q$ is also a CWC$(n,d,w)_q$.

Johnson-type bounds for constant-weight codes and constant-composition codes have been derived.

\begin{lemma}[Svanstr\"{o}m \cite{Svanstrom:1999a}, Svanstr{\"o}m \etal{} \cite{Svanstrometal:2002}]\label{bound-cwcccc}
\begin{align*}
A_q(n,d,\w)&\le \floor{\frac{n}{w_1}A_q(n-1,d,[w_1-1,\ldots,w_{q-1}])},\\
A_q(n,d,w)& \leq\floor{\frac{(q-1)n}{w}A_q(n-1,d,w-1)}.\\
\end{align*}
\end{lemma}

%

Apply the fact that $A_q(n,2w,\w)= \floor{n/w}$ and $A_q(n,2w,w) =\floor{n/w}$
(see Fu \etal{} \cite{Fuetal:1998}, Chee \etal{} \cite{Cheeetal:2008b})
 to Lemma \ref{bound-cwcccc}, 
we have the following upper bounds:
\begin{align}
A_q(n,2w-2,\w)&\le \floor{\frac{n}{w_1}\floor{\frac{n-1}{w-1}}},\label{ccc1}\\
A_q(n,2w-3,\w)& \le 
\begin{cases}
  \floor{\frac{n}{w_1}\floor{\frac{n-1}{w_1-1}}}, & \mbox{if $w_1>w_2$,}\\[6pt]
   \floor{\frac{n}{w_1}\floor{\frac{n-1}{w_1}}}, & \mbox{otherwise,}
   \end{cases}\label{ccc2}\\
A_q(n,2w-2,w)&\le   \floor{\frac{(q-1)n}{w}\floor{\frac{n-1}{w-1}}},\label{cwc1}\\
A_q(n,2w-3,w)&\le \floor{\frac{(q-1)n}{w}\floor{\frac{(q-1)(n-1)}{w-1}}}.\label{cwc2}
\end{align}

In this paper, we show that the above inequalities are exact provided that $n$ is sufficiently large
and $n$ satisfies certain congruence conditions. 
To do so, we apply the theory of decompositions of edge-colored digraphs to construct optimal codes
meeting the above bounds.

We list previous similar asymptotic or exact results.

\begin{enumerate}[(i)]
\item Results for $A_q(n,d,\w)$ are known

{\small
\begin{enumerate}[(a)]
\item for all $\w$ and $d=2w-1$ \cite{Cheeetal:2010a};
\item for all $d$ where $w\le 3$ \cite{Svanstrom :2000,Cheeetal:2010a};
\item for $(q,d,w)=(3,5,4)$ \cite{GaoGe:2011}.
\end{enumerate}
}
\item Results for $A_q(n,d,w)$ are known 
{\small
\begin{enumerate}[(a)]
\item for all $w$ and $d=2w-1$ \cite{CheeLing:2007,Cheeetal:2010a};
\item for all $q$ and $(d,w)\in\{(3,2),(4,3),(5,3)\}$ \cite{CheeLing:2007,Cheeetal:2008,Cheeetal:2013htp};
\item for $(q,d,w)\in\{(3,5,4),(3,6,4),(4,5,4)\}$ \cite{ZhangGe:2010,Zhangetal:2012,ZhangGe:2013}. 
\end{enumerate}
}
\end{enumerate}

\subsection{Multiply Constant-Weight Codes}

Consider a binary code $\C\subseteq \bz_2^{[m]\times [n]}$ of constant-weight $mw$ and distance $d$.
The code $\C$ is said to be of {\em multiply constant-weight $w$} if for $\su\in \C$, $i\in [m]$,
the subword $(\su_{i,j})_{j\in [n]}$ is of constant-weight $w$.
Denote such a code by MCWC$(m,n,d,w)$.
Similarly, the maximum size of a MCWC$(m,n,d,w)$ is given by $M(m,n,d,w)$ and 
a multiply constant-weight code attaining the maximum size is said to be {\em optimal}.
A Johnson-type bound can again be derived.

\begin{lemma}[Chee \etal{} \cite{Cheeetal:preprint2013}]
\begin{equation}
M(m,n,2mw-2,w)  \le \floor{\frac{n}{w}\floor{\frac{n}{w}}}.\label{mcwc}
\end{equation}
\end{lemma}

Again, we verify \eqref{mcwc} is exact provided $n$ is sufficiently large
 and satisfies certain congruence conditions. 
Next, we describe the main tool in our constructions of optimal codes.

\subsection{Decomposition of Edge-Colored Complete Digraphs}

Denote the set of all ordered pairs of a finite set $X$ with distinct components by $\overline{\binom{X}{2}}$.
An {\em edge-colored digraph} is a triple $G=(V,C,E)$, where $V$ is a finite set of {\em vertices}, $C$ is a finite set of {\em colors}
and $E$ is a subset of $\overline{\binom{V}{2}}\times C$. Members of $E$ are called {\em edges}. The {\em complete edge-colored digraph} on $n$ vertices with $r$ colors, denoted by $K^{(r)}_n$, is the edge-colored digraph $(V,C,E)$, where $|V|=n$, $|C|=r$ and $E=\overline{\binom{V}{2}}\times C$.

A family $\F$ of edge-colored subgraphs of an edge-colored digraph $K$ is a {\em decomposition} of $K$ if every edge of $K$ belongs to exactly one member of $\F$. Given a family of edge-colored digraphs $\G$, a decomposition $\F$ of $K$ is a {\em $\G$-decomposition of $K$} if each edge-colored digraph in $\F$ is isomorphic to some $G\in \G$.
Furthermore, a $\G$-decomposition of $K$ is said to be {\em superpure} 
if any two distinct edge-colored subgraphs in $\F$ share at most two vertices.

Lamken and Wilson \cite{LamkenWilson:2000} exhibited the asymptotic existence of 
decompositions of $K^{(r)}_n$ for a fixed family of digraphs.
Hartmann \cite{Hartmann:2002} later extended their results to superpure decompositions.
To state the theorems, we require more concepts.


Consider an edge-colored digraph $G=(V,C,E)$ with $|C|=r$. Let $((u,v),c)\in E$ denote a directed edge from $u$ to $v$, colored by $c$. For any vertex $u$ and color $c$, define the {\em indegree} and {\em outdegree} of $u$ with respect to $c$, to 
be the number of directed edges of color $c$ entering and leaving $u$ respectively.
Then for vertex $u$, we define the {\em degree vector} of $u$ in $G$, denoted by $\tau(u,G)$, to be the vector of length $2r$,
$\tau(u,G) \triangleq(\ind_1(u,G),\outd_1(u,G),\dots,\ind_r(u,G),\outd_r(u,G))$.
Define $\alpha(\G)$ to be the greatest common divisor of the integers $t$ such that the $2r$-vector $(t,t,\ldots,t)$ 
is a nonnegative integral linear combination of the degree vectors $\tau(u,G)$ as $u$ ranges over all vertices of all digraphs  $G\in \G$.

For each $G=(V,C,E)\in \G$, let $\mu(G)$ be the {\em edge vector} of length $r$ given by 
$\mu(G) \triangleq(m_1(G),m_2(G),\dots,m_r(G))$ where 
$m_i(G)$ is the number of edges with color $i$ in $G$. 
We denote by $\beta(\G)$ the greatest common divisor of the integers $m$ such that $(m,m,\dots,m)$ is a nonnegative integral linear combination of the vectors $\mu(G)$, $G\in\G$.
Then $\G$ is said to be {\em admissible} if $(1,1,\ldots,1)$ can be expressed as a positive rational combination of  the vectors $\mu(G)$, $G\in\G$.



Below is the main theorem we allude to in our proofs.

\begin{theorem} [Lamken and Wilson \cite{LamkenWilson:2000}, Hartmann \cite{Hartmann:2002}]
\label{main} Let $\G$ be an admissible family of edge-colored digraphs with $r$ colors. 
Then there exists a constant $n_0$ (resp. $n_1$) 
such that a (resp. superpure) $\G$-decomposition of $K_n^{(r)}$ exists for every 
$n\geq n_0$ (resp. $n\ge n_1$) satisfying:
$n(n-1)\equiv 0\pmod{\beta(\G)}$ and
$n-1\equiv 0\pmod{\alpha(\G)}$.

\end{theorem}

We apply Theorem \ref{main} to construct codes that meet the upper bounds \eqref{ccc1}--\eqref{mcwc}.
To do so, we have two main steps.
\begin{enumerate}[(A)] 
\item 
Define a family $\G$ of edge-colored digraphs on a set of $r$ colors. 
We then show that a $\G$-decomposition of $K_n^{(r)}$ results in a code of length $n$ satisfying certain
weight and distance properties.
\item Compute $\alpha(\G)$ and $\beta(\G)$ and hence determine the congruences classes 
that $n$ needs to satisfy.
\end{enumerate}

In this paper, we focus on step (A), describing the construction of $\G$ and 
establishing the correspondence to certain codes. 
As the computations in step (B) are usually tedious and analogous, 
we exhibit one computation for illustrative purposes and 
defer other computations to the 
appendices. 

\section{Asymptotically Exact Johnson-Type Bounds for Constant-Composition Codes}
\label{sec:ccc}

Fix $\w=[w_1,w_2,\dots,w_{q-1}]$ with $w_1\geq \dots \geq w_{q-1}>0$ and let $w=\sum_{i=1}^{q-1}w_i$. 
In this section, we construct infinite families of optimal CCC$(n,d,\w)_q$ for $d=2w-2$ or $d=2w-3$,
and establish the asymptotic exactness of \eqref{ccc1} and \eqref{ccc2}.

\subsection{When Distance $d=2w-2$}

Consider the following characterization
of codes of constant-weight $w$ with distance $2w-2$.
Note that since the constant-composition codes are constant-weight codes, 
the lemma is applicable for both classes of codes.

\begin{lemma}\label{ccc2w-2} 
The following are necessary and sufficient for a code $\C$ 
of constant-weight $w$ to have distance $2w-2$:
\begin{enumerate}[{\rm (C1)}]
\item For $i\in [q-1]$, the ordered pairs in the set 
$\{(x,y):\su_{x}=i,y\in \supp(\su)\setminus\{x\}, \su\in \C\}$
are distinct.
\item For any $\su,\sv\in \C$, $|\supp(\su)\cap \supp(\sv)|\leq 2$. 
\end{enumerate}
\end{lemma}

%

\noindent{\bf Definition of the family $\G(\w)$}.
For fixed $\w$, define an edge-colored digraph $G(\w)=(V[\w],C[\w],E[\w])$, where

{\small
\begin{align*}
V[\w]&\triangleq\{x_{ij}:i\in[q-1],j\in [w_i]\}; \\
C[\w]&\triangleq[q-1]; \\
E[\w]&\triangleq\left\{((x_{ij},x_{ij'}),i):i\in[q-1],(j,j)'\in\overline{\binom{[w_i]}{2}}\right\}\cup \\
&~\left\{((x_{ij},x_{i'j'}),i):(i,i')\in\overline{\binom{[q-1]}{2}},j\in[w_i],j'\in[w_{i'}]\right\}. 
\end{align*}
}
Let $s$ be the largest integer such that $w_1=w_2=\dots=w_s$. 
For each $s+1\le i\le q-1$, let $G_i$ be an edge-colored digraph with two vertices 
$y_i$, $z_i$ and one directed edge with color $i$ from $y_i$ to $z_i$. 
Then $\G(\w)=\{G(\w)\}\cup \{G_i:  s+1\le i\le q-1\}$.


\begin{example}
Let $\w=[3,2]$. The edge-colored digraph $G(\w)$ is given below, 
where the solid lines denote directed edges with color $1$, the dotted lines denote the directed edges with color $2$, and ``$\xymatrix {\ar@{<->}[r] &}$'' (``$\xymatrix {\ar@{<.>}[r] &}$'') denotes the two directed edges with the same color with one in each direction.
\end{example}

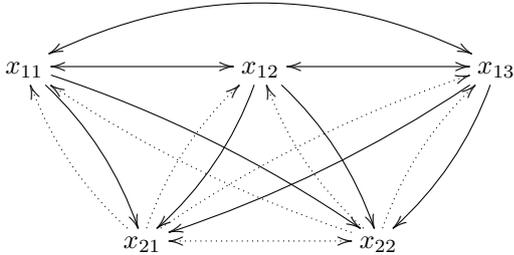
\begin{equation*}
 \xymatrix{
x_{11} \ar@{<->}[rr] \ar@/^2pc/@{<->}[rrrr] \ar@/^/[ddr] \ar@/^/[ddrrr] & & 
x_{12} \ar@{<->}[rr] \ar@/^/[ddl] \ar@/^/[ddr] & & 
x_{13} \ar@/^/[ddlll] \ar@/^/[ddl] \\ \\
 & x_{21} \ar@/^/@{.>}[uul] \ar@/^/@{.>}[uur] \ar@/^/@{.>}[uurrr] \ar@{<.>}[rr] & & 
 x_{22} \ar@/^/@{.>}[uulll] \ar@/^/@{.>}[uul] \ar@/^/@{.>}[uur]
}
\end{equation*}
Then $G_2$ is the digraph $\xymatrix {y_2\ar@{.>}[r] &z_2}$ and 
the family of digraphs is given by $\G(\w)=\{G(\w),G_2\}$.

\vspace{2mm}

\noindent{\bf Construction of an optimal CCC$(n,2w-2,\w)_q$}. 
Suppose a superpure $\G(\w)$-decomposition of $K_n^{(q-1)}$ exists. 
For each $F$ isomorphic to $G(\w)$, there is a unique partition $V(F)=\bigcup_{i=1}^{q-1}S_i$ 
so that the edges from $x$ in $F$ has color $i$ if and only if $x\in S_i$. 
Then construct one codeword $\su$ with support $V(F)$ such that $\su_x=i$ for $x\in S_i$. 
Hence, $\su$ has composition $\w$.

Since we have a $\G(\w)$-decomposition of $K^{(q-1)}_n$, then (C1) of Lemma~\ref{ccc2w-2} is satisfied.
Furthermore, since the decomposition is superpure, (C2) is also met. 
Hence, we obtain a CCC$(n,2w-2,\w)_q$ code of size ${n(n-1)}/(w_1(w-1))$. 
The code meets the upper bound in \eqref{ccc1} and is therefore optimal.

\begin{example} Let $\w=[2,1]$. Consider the following superpure $\G(\w)$-decomposition of $K^{(2)}_5$:

{
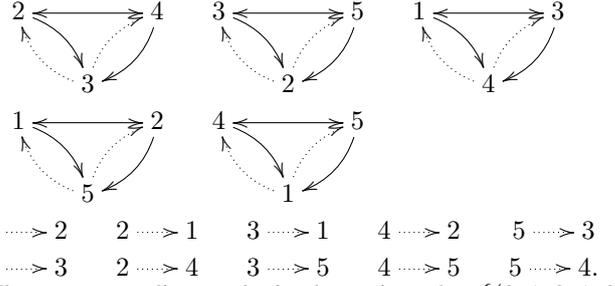
\begin{tabular}{ccc}
 \xymatrix@=1.5em{
2 \ar@{<->}[rr] \ar@/^/@{->}[dr]& & 
4 \ar@/^/@{->}[dl]  \\
 & 3\ar@/^/@{.>}[ul]\ar@/^/@{.>}[ur]
 }
&
 \xymatrix @= 1.5em{
3 \ar@{<->}[rr] \ar@/^/@{->}[dr]& & 
5 \ar@/^/@{->}[dl]  \\
 & 2\ar@/^/@{.>}[ul]\ar@/^/@{.>}[ur]
 }
&
 \xymatrix@= 1.5em{
1 \ar@{<->}[rr] \ar@/^/@{->}[dr]& & 
3 \ar@/^/@{->}[dl]  \\
 & 4\ar@/^/@{.>}[ul]\ar@/^/@{.>}[ur]
 }
\\
\xymatrix@= 1.5em{
1 \ar@{<->}[rr] \ar@/^/@{->}[dr]& & 
2 \ar@/^/@{->}[dl]  \\
 & 5\ar@/^/@{.>}[ul]\ar@/^/@{.>}[ur]
 }
&
 \xymatrix@= 1.5em{
4 \ar@{<->}[rr] \ar@/^/@{->}[dr]& & 
5 \ar@/^/@{->}[dl]  \\
 & 1\ar@/^/@{.>}[ul]\ar@/^/@{.>}[ur]
 }
\end{tabular}
\begin{tabular}{ccccc}
 \xymatrix@=1.5em{
1 \ar@{.>}[r]& 2
 }
&
\xymatrix@=1.5em{
2 \ar@{.>}[r]&1
 }
&
\xymatrix@=1.5em{
3 \ar@{.>}[r]&1
 }
&
\xymatrix@=1.5em{
4 \ar@{.>}[r]&2
 }
&
\xymatrix@=1.5em{
5 \ar@{.>}[r]&3
 }
 \\
  \xymatrix@=1.5em{
1 \ar@{.>}[r]& 3
 }
&
\xymatrix@=1.5em{
2 \ar@{.>}[r]&4
 }
&
\xymatrix@=1.5em{
3 \ar@{.>}[r]&5
 }
&
\xymatrix@=1.5em{
4 \ar@{.>}[r]&5
 }
&
\xymatrix@=1.5em{
5 \ar@{.>}[r]&4.
 }
\end{tabular}
}

The corresponding code is then given by 
$\{(0, 1, 2, 1, 0)$, $(0, 2, 1, 0, 1), (1, 0, 1, 2, 0), (1, 1, 0, 0, 2),(2, 0, 0, 1, 1)\}$, 
which is indeed a CCC$(5,4,[2,1])_3$.
\end{example}

\noindent{\bf Computation of $\alpha(\G(\w))$ and $\beta(\G(\w))$}.
First consider the digraph $G(\w)$. 
Observe that for $i\in[q-1], k\in[w_i]$ we have $\ind_i(x_{ik},G(\w))=w_i-1$, $\outd_i(x_{ik},G(\w))=w-1$,
$\ind_j(x_{ik},G(\w))=w_j$ and $\outd_j(x_{ik},G(\w))=0$ for $j\ne i$.
Consider $G_i$ for $s+1\le i\le q$.
Then $\ind_i(z_i,G_i)=\outd_i(y_i,G_i)=1$
and all other indegrees and outdegrees are zero.

Let $a=\gcd(w_1,w)$. Pick $t=\floor{w/w_1}$ so that $0\le w-tw_1<w_1$.
Observe also that $t\ge s$.
Consider the vector
\begin{equation*}
\upsilon=\frac{w-tw_1}{a}\tau(x_{(t+1)1},G(\w))+\frac{w_1}{a}\sum_{i=1}^t\tau(x_{i1},G(\w)).
\end{equation*}
For $j\in [q-1]$, let $\ind_j(\upsilon)$ and $\outd_j(\upsilon)$ denote the coordinates in $\upsilon$
corresponding to the summation of the indegrees or outdegrees with respect to color $j$.
Then we have 
\begin{align*}
\ind_j(\upsilon) &=
\begin{cases}
\frac{ww_j-w_1}{a}=\frac{w_1(w-1)}{a}, & \text{for } j\in [s],\\
\frac{ww_j-w_1}{a}<\frac{w_1(w-1)}{a}, & \text{for } s+1\le j\le t,\\
 \frac{ww_{t+1}-w+tw_1}{a}< \frac{w_{t+1}(w-1)}{a}, & \text{for } j=t+1,\\
\frac{ww_{j}}{a}, & \text{otherwise},
\end{cases}
\\
\outd_j(\upsilon) &=
\begin{cases}
\frac{w_1(w-1)}{a}, & \text{for } j\in [t],\\
 \frac{(w-tw_1)(w-1)}{a}< \frac{w_{1}(w-1)}{a}, & \text{for } j=t+1,\\
 0, & \text{otherwise}.
\end{cases}
\end{align*}
Observe that the first $2s$ coordinates of $\upsilon$ are $w_1(w-1)/a$ and 
all other coordinates have values at most $w_1(w-1)/a$.
Adding to $\upsilon$ a suitable nonnegative integral combinations of $\tau(y_i,G_i)$'s and $\tau(z_i,G_i)$'s,
we get $w_1(w-1)/a(1,1,\ldots,1)$. 
Hence, we conclude that $\alpha(\G(\w))=w_1(w-1)/a$.

Next, consider the edge vector $\mu(G)$ with $G\in\G(\w)$.
For $G(\w)$, $m_i(G(\w))=w_i(w-1)$ for $i\in [q-1]$, while 
for $G_i$ with $s+1\le i\le q-1$, $m_i(G_i)=1$ and $m_j(G_i)=0$ for $j\ne i$. 
Hence, $\beta(\G(\w))=w_1(w-1)$.

Applying Theorem \ref{main}, we obtain our first asymptotic result.

\begin{proposition}\label{main:ccc1} 
Fix $\w$ and let $w=\sum_{i=1}^{q-1}w_i$. There exists an integer $n_0$ such that 
\begin{equation*}
A_q(n,2w-2,\w)=\frac{n(n-1)}{w_1(w-1)}
\end{equation*} 
\noindent for all $n\geq n_0$ satisfying 
$n(n-1)\equiv  0\pmod{w_1(w-1)}$ and 
$n-1\equiv  0\pmod{w_1(w-1)/a}$, 
where $a=\gcd(w_1,w)$.
\end{proposition}

\subsection{When Distance $d=2w-3$}

We have the following analogous characterization of codes of constant-weight $w$
with distance $2w-3$.
\begin{lemma}
\label{ccc2w-3} The following are sufficient for a code $\C$ of weight $w$ to have distance $2w-3$.
\begin{enumerate}
\item[{\rm (C3)}] For $i,j\in [q-1]$, the ordered pairs in the set
$\{(x,y):\su_x=i,\su_y=j, x\ne y,\su\in \C\}$
are distinct.
\item[{\rm (C4)}] For any $\su,\sv\in \C$, $|\supp(\su)\cap \supp(\sv)|\leq 2$. 
\end{enumerate}
\end{lemma}

\noindent{\bf Definition of $\G^*(\w)$}. 
For fixed $\w$, define an edge-colored digraph $G^*(\w)=(V(\w),C(\w),E(\w))$, where

{\footnotesize
\begin{align*}
V(\w)&\triangleq\{x_{ij}:i\in[q-1],j\in [w_i]\}; \\
C(\w)&\triangleq[q-1]\times [q-1]; \\
E(\w)&\triangleq\left\{((x_{ij},x_{ij'}),(i,i)):i\in[q-1],(j,j)'\in\overline{\binom{[w_i]}{2}}\right\}\cup \\
&\left\{((x_{ij},x_{i'j'}),(i,i')):(i,i')\in\overline{\binom{[q-1]}{2}},j\in[w_i],j'\in[w_{i'}]\right\}. 
\end{align*}
}

For $i,j\in[q-1]$, let $G^*_{ij}$ be 
a digraph with vertices
$y_i$, $z_j$ and one directed edge with color $(i,j)$ from $y_i$ to $z_j$.
To define $\G^*(\w)$, we have two cases depending on whether $w_1=w_2$:
\begin{enumerate}[(i)]
\item When $w_1> w_2$, let $r$ be the largest integer such that $w_2=\dots=w_r=w_1-1$.
Then set $\G^*(\w)=\{G^*(\w)\}\cup \{G^*_{ij}:(i,j)\in [q-1]\times[q-1]\setminus \{(1,1),(1,2),(1,3),\dots,(1,r),(2,1),(3,1)\dots,(r,1)\}$.
\item When $w_1=w_2$, let $r$ be the largest integer such that $w_1=\dots=w_r$. 
Then set $\G^*(\w)=\{G^*(\w)\}\cup \left\{G^*_{ij}:(i,j)\in [q-1]\times [q-1]\setminus \overline{\binom{[r]}{2}}\right\}$.
\end{enumerate} 

\begin{example}
Let $\w=[3,2]$. The edge-colored digraph $G^*(\w)$ is given below, 
where the lines ``$\xymatrix {\ar[r] &}$'', ``$\xymatrix {\ar@{-->}[r] &}$'', ``$\xymatrix {\ar@{.>}[r] &}$'' and ``$\xymatrix {\ar@{~>}[r] &}$'' 
denote directed edges with color $(1,1)$, $(1,2)$, $(2,1)$ and $(2,2)$ respectively. 

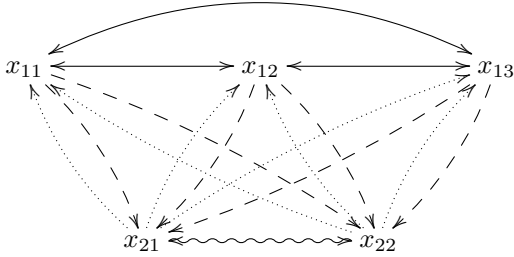
\begin{equation*}
 \xymatrix{
x_{11} \ar@{<->}[rr] \ar@/^2pc/@{<->}[rrrr] \ar@/^/@{-->}[ddr] \ar@/^/@{-->}[ddrrr] & & 
x_{12} \ar@{<->}[rr] \ar@/^/@{-->}[ddl] \ar@/^/@{-->}[ddr] & & 
x_{13} \ar@/^/@{-->}[ddlll] \ar@/^/@{-->}[ddl] \\ \\
 & x_{21} \ar@/^/@{.>}[uul] \ar@/^/@{.>}[uur] \ar@/^/@{.>}[uurrr] \ar@{<~>}[rr] & & 
 x_{22} \ar@/^/@{.>}[uulll] \ar@/^/@{.>}[uul] \ar@/^/@{.>}[uur] 
}
\end{equation*}
Then $G^*_{22}$ is the digraph $\xymatrix {y_2\ar@{~>}[r] &z_2}$ and 
the family of digraphs is given by 
$\G^*(\w)=\{G^*(\w),G^*_{22}\}$.
\end{example}
\vspace{2mm}

\noindent{\bf Construction of an optimal CCC$(n,2w-3,\w)_q$}. 
Suppose a superpure $\G^*(\w)$-decomposition of $K_n^{(q-1)^2}$ exists. 
For $F$ isomorphic to $G^*(\w)$, there is a unique partition $V(F)=\bigcup_{i=1}^{q-1}S_i$ 
so that the edges from $x$ to $y$ in $F$ has color $(i,j)$ 
if $x\in S_i$ and $y\in S_j$. 
Construct a codeword $\su$ with support $V(F)$ such that $\su_x=i$ for $x\in S_i$. 
So, $\su$ has composition $\w$.

Since we have a superpure $\G^*(\w)$-decomposition of $K^{(q-1)^2}_n$, (C3) and (C4) of Lemma~\ref{ccc2w-3} are satisfied.
Hence, we obtain a CCC$(n,2w-3,\w)_q$ code. 
Let $M$ be the number of digraphs isomorphic to $G^*(\w)$. It is easy to see that 
$M=n(n-1)/(w_1(w_1-1))$ if $w_1\ne w_2$ and $M=n(n-1)/w_1^2$, otherwise.
The code is optimal by the upper bound \eqref{ccc2}.

We compute $\alpha(\G^*(\w))$ and $\beta(\G^*(\w))$ in the full paper.
\begin{align*}
\alpha(\G^*(\w))&=
\begin{cases}
w_1(w_1-1), & \text{if }w_1> w_2,\\
w_1, & \text{otherwise}.
\end{cases}\\
\beta(\G^*(\w))&=
\begin{cases}
w_1(w_1-1), & \text{if } w_1> w_2,\\
w_1^2, & \text{otherwise}.
\end{cases}
\end{align*}

Applying Theorem \ref{main} we obtain the following proposition.

\begin{proposition}\label{main:ccc2} 
Fix $\w$ and let $w=\sum_{i=1}^{q-1}w_i$. There exists an integer $n_0$ such that 
\begin{equation*}
A_q(n,2w-3,\w)=
\begin{cases}
\frac{n(n-1)}{w_1(w_1-1)}, & \text{if }w_1> w_2,\\
\frac{n(n-1)}{w_1^2}, & \text{otherwise},
\end{cases}
\end{equation*} 
for all $n\geq n_0$ satisfying 
\begin{enumerate}[{\rm (i)}]
\item $n-1\equiv  0\pmod{w_1(w_1-1)}$, if $w_1> w_2$, 
\item $n-1\equiv  0\pmod{w_1^2}$, otherwise. 
\end{enumerate}
\end{proposition}

\section{Asymptotically Exact Johnson-Type Bounds for $q$-ary Constant-Weight Codes}
\label{sec:cwc}

In this section, we construct infinite families of optimal CWC$(n,d,w)_q$ codes for $d=2w-2$ or $d=2w-3$,
and establish the asymptotic exactness of \eqref{cwc1} and \eqref{cwc2}.

Note that constant-composition codes are special instances of constant-weight codes.
Hence, we make use of the digraphs defined in Section \ref{sec:ccc} to form
the desired families for constant-weight codes.
Specifically, consider the set of compositions,

{\footnotesize
\begin{equation*}
W=\left\{[w_1,w_2,\ldots,w_{q-1}]: 0\le w_i\le w \text{ for } i\in[q-1], \sum_{i=1}^{q-1}w_i=w\right\}.
\end{equation*}
}
Note here we do not require all values in the composition to be positive and 
the composition to be monotone decreasing. 
Furthermore, we can extend the constructions in Section \ref{sec:ccc}
to define $G(\w)$ and $G^*(\w)$ for all $\w\in W$.
Define the following families of digraphs 
\begin{equation*}
\G(w)\triangleq\bigcup_{\w\in W} G(\w)\text{ and }\G^*(w)\triangleq\bigcup_{\w\in W} G^*(\w).
\end{equation*}

Similar to constructions in Section \ref{sec:ccc}, we have that
a superpure $\G(w)$-decomposition of $K^{(q-1)}_n$
and a superpure $\G^*(w)$-decomposition of $K^{(q-1)^2}_n$
yield a CWC$(n,2w-2,w)_q$ and a CWC$(n,2w-3,w)_q$ respectively.

The number of digraphs in a superpure  $\G(w)$-decomposition of $K^{(q-1)}_n$
 is $(q-1)n(n-1)/(w(w-1))$
since the total number of directed edges is $(q-1)n(n-1)$. 
In the full paper 
we used the methods of Lamken and Wilson to
establish the following.

\begin{proposition}\label{main:cwc1} 
Fix $w$. 
There exists an integer $n_0$ such that 
\begin{equation*}
A_q(n,2w-2,w)=\frac{(q-1)n(n-1)}{w(w-1)}
\end{equation*}
for all $n\ge n_0$ satisfying $(q-1)n(n-1)\equiv 0\pmod{w(w-1)}$ and $n-1\equiv 0 \pmod{w-1}$.
\end{proposition}

On the other hand, $\G^*(w)$ corresponds to the family of digraphs constructed in the proof of 
 \cite[Theorem 8.1]{LamkenWilson:2000}. 
 Therefore, we have the following proposition.

\begin{proposition}\label{main:cwc2} 
Fix $w$. 
There exists an integer $n_0$ such that 
\begin{equation*}
A_q(n,2w-3,w)=\frac{(q-1)^2n(n-1)}{w(w-1)}
\end{equation*}
for all $n\ge n_0$ satisfying $(q-1)^2n(n-1)\equiv 0\pmod{w(w-1)}$ and $(q-1)(n-1)\equiv 0 \pmod{w-1}$.
\end{proposition}

\section{Asymptotically Exact Johnson-Type Bounds for Multiply Constant-Weight Codes}
\label{sec:mcwc}

In this section, we construct an infinite family of optimal MCWC$(m,n,2mw-2,w)$ 
and establish the asymptotic exactness of \eqref{mcwc}.
Fix $m$ and $w$. Using digraphs from Section \ref{sec:ccc}, we define the edge-colored digraph 
$H^*(m,w)\triangleq G^*(\w)$, with $w_1=w_2=\cdots=w_{m}=w$ (here $q-1=m$).
Define $\HH^*(m,w)=\{H^*(m,w)\}\cup\{G^*_{ii}:i\in [m]\}$. 
\vspace{2mm}

\noindent{\bf Construction of an optimal MCWC$(m,n,2mw-2,w)$}.
Suppose an $\HH^*(m,w)$-decomposition of  $K^{(m^2)}_n$ exists.
 For $F$ isomorphic to $H^*(m,w)$, there is a unique partition $V(F)=\bigcup_{i=1}^{m}S_i$ 
so that the edges from $x$ to $y$ in $F$ has color $(i,j)$ if and only if $x\in S_i$ and $y\in S_j$. 
Then construct one codeword $\su$ such that $\su_{(i,x)}=1$ for $i\in [m]$ and $x\in S_i$.

Since we have an $\HH^*(m,w)$-decomposition of $K^{(m^2)}_n$, 
the distance of the code is $2mw-2$.
Hence, we have an MCWC$(m,n,2mw-2,w)$ of size $n(n-1)/w^2$
that meets the upper bound given by \eqref{mcwc}.
Similar computations yield $\alpha(\HH^*(m,w))=w$ and $\beta(\HH^*(m,w))=w^2$.
Applying Theorem \ref{main} we obtain the following.

\begin{proposition}\label{main:mcwc} 
Fix $m$ and $w$. 
There exists an integer $n_0$ so that 
\begin{equation*}
M(m,n,2mw-2,w)=\frac{n(n-1)}{w^2}
\end{equation*}
for all $n\ge n_0$ satisfying $n-1\equiv 0 \pmod{w^2}$.
\end{proposition}

\section{Conclusion}

We verify that Johnson-type bounds are asymptotically exact 
for several generalizations of binary constant-weight codes.
This was achieved via an interesting application of 
superpure decompositions of edge-colored digraphs.

Observe that in Propositions \ref{main:ccc1}, \ref{main:ccc2},  \ref{main:cwc1}, \ref{main:cwc2} and \ref{main:mcwc} 
Johnson-type bounds are shown to be exact for sufficiently large $n$
satisfying certain congruence classes. 
We hypothesize that the Johnson-type bounds are tight up to an additive constant
for sufficiently large $n$. Specifically, we make the following conjecture.

\vspace{3mm}

\noindent{\bf Conjecture.} Fix $q$, $\w$, $w$, $m$ and let $d\in\{2w-2,2w-3\}$.
Define $U_q(n,d,\w)$, $U_q(n,d,w)$, $U(m,n,d,w)$ to be the upper bounds given by
\eqref{ccc1} -- \eqref{mcwc}. Then
\begin{align*}
A_q(n,d,\w) &= U_q(n,d,\w)-O(1),\\
A_q(n,d,w) &= U_q(n,d,w)-O(1),\\
M(m,n,2mw-2,w) &= U(m,n,2mw-2,w)-O(1).
\end{align*}
The case for binary constant-weight codes, that is, $q=2$ and $d=2w-2$,
has been verified by Chee \etal{} \cite{Cheeetal:2013pairs}.

\appendices

\section{On $\G^*(\w)$ and the Proof of Proposition \ref{main:ccc2}}\label{compute:ccc2}


We compute $\alpha(\G^*(\w))$ and $\beta(\G^*(\w))$ defined in section \ref{sec:ccc}.

\subsection{When $ w_1> w_2$.}

Recall $r$ is the largest integer such that $w_2=\dots=w_r=w_1-1$
and $\G^*(\w)=\{G^*(\w)\}\cup \{G^*_{ij}:(i,j)\in [q-1]\times[q-1]\setminus \{(1,1),(1,2),(1,3),\dots,(1,r),(2,1),(3,1)\dots,(r,1)\}$.

Consider the digraph $G^*(\w)$. Observe that for $i\in [q-1]$ and $k\in [w_i]$,
we have $\ind_{(i,i)}(x_{ik}, G^*(\w))=\outd_{(i,i)}(x_{ik}, G^*(\w))=w_i-1$, 
$\ind_{(j,i)}(x_{ik}, G^*(\w))=\outd_{(i,j)}(x_{ik}, G^*(\w))=w_j$ for $j\ne i$ and 
and the other indegrees and outdegrees are zero. 
For the digraph $G_{ij}$, we have $\outd_{(i,j)}(y_{ij},G_{ij})=\ind_{(j,i)}(z_{ij},G_{ij})=1$
and the other indegrees and outdegrees are zero.

Consider the vector
\begin{equation*}
\upsilon=w_1\tau(x_{11},G^*(\w))+(w_1-1)\sum_{i=2}^{[q-1]}\tau(x_{i1},G^*(\w)).
\end{equation*}
For $(i,j)\in [q-1]\times [q-1]$, let $\ind_{(i,j)}(\upsilon)$ and $\outd_{(i,j)}(\upsilon)$ denote the coordinates in $\upsilon$
corresponding to the summation of the indegrees and outdegrees with respect to color $(i,j)$.
Then we have 
\begin{align*}
\ind_{(i,j)}(\upsilon) &=
\begin{cases}
w_1(w_1-1), & \text{for } i=j=1,\\
w_1(w_1-1), & \text{for } i=1,j\ne 1,\\
w_1w_i & \text{for } i\ne 1,j=1,\\
(w_1-1)(w_i-1) & \text{for } i=j, i\ne 1,\\
w_i(w_1-1), & \text{otherwise},
\end{cases}
\\
\outd_{(i,j)}(\upsilon) &=
\begin{cases}
w_1(w_1-1), & \text{for } i=j=1,\\
w_1w_j, & \text{for } i=1,j\ne 1,\\
w_1(w_1-1) & \text{for } i\ne 1,j=1,\\
(w_1-1)(w_i-1) & \text{for } i=j, i\ne 1,\\
w_j(w_1-1), & \text{otherwise},
\end{cases}
\end{align*}

Observe that the coordinates of $\upsilon$ corresponding to 
indegrees and outdegrees with respect to colors in 
 $\{(1,1),(1,2),(1,3),\dots,(1,r),(2,1),(3,1)\dots,(r,1)\}$
 have value $w_1(w_1-1)$.
All other coordinates have values at most $w_1(w-1)$.
Adding to $\upsilon$ a suitable nonnegative integral combinations of $\tau(y_{ij},G_{ij})$'s and $\tau(z_{ij},G_{ij})$'s,
we obtain $(w_1(w-1),w_1(w-1),\ldots,w_1(w-1))$.
Hence, we conclude that $\alpha(\G)=w_1(w-1)$.

Next, consider the edge vector $\mu(G)$ with $G\in\G^*(\w)$.
For $G^*(\w)$, we have $m_{(i,i)}(G^*(\w))=w_i(w_i-1)$ and $m_{(i,j)}(G^*(\w))=w_iw_j$ for $j\ne i$.
On the other hand, for $G_{ij}$  we have $m_{(i,j)}(G_{ij})=1$ and $m_{(i',j')}(G_{ij})=0$ for $(i',j')\ne (i,j)$. 
Hence, $\beta(\G^*(\w))=w_1(w-1)$.

\subsection{When $w_1=w_2$.}

Recall $r$ is the largest integer such that $w_1=\dots=w_r$ and
 $\G^*(\w)=\{G^*(\w)\}\cup \left\{G_{ij}:(i,j)\in [q-1]\times [q-1]\setminus \overline{\binom{[r]}{2}}\right\}$. 
 
 Here, we consider the vector
 $\upsilon'=\sum_{i=1}^{q-1} \tau(x_{i1},G^*(\w))$.
 Then the coordinates of $\upsilon'$ corresponding to 
indegrees and outdegrees with respect to colors in 
 $\overline{\binom{[r]}{2}}$
 have value $w_1$.
All other coordinates have values at most $w_1$.
Adding to $\upsilon$ a suitable nonnegative integral combinations of $\tau(y_{ij},G_{ij})$'s and $\tau(z_{ij},G_{ij})$'s,
we conclude that $\alpha(\G^*(\w))=w_1$.
Similar to the case where $w_1\ne w_2$, we have $\beta(\G^*(\w))=w_1^2$.

\section{On $\G(w)$ and the Proof of Proposition \ref{main:cwc1}}\label{compute:cwc1}

Unfortunately, it is not straightforward to determine 
$\alpha(\G(w))$ and $\beta(\G(w))$.
Instead, we follow the methodology  of Lamken and Wilson (see \cite[Theorem 8.1]{LamkenWilson:2000})
to complete the proof of Proposition \ref{main:cwc1}.
Specifically, we show that provided $(q-1)n(n-1)\equiv 0\pmod{w(w-1)}$
$n-1\equiv 0\pmod{w-1}$, we have:
\begin{enumerate}[(i)]
\item the vector $n(n-1)(1,1,\ldots,1)$ is 
an integral linear combination of vectors in $\{\mu(G): G\in\G(w)\}$;
\item the vector $(n-1)(1,1,\ldots,1)$ is 
an integral linear combination of vectors in $\{\tau(u,G): u\in G,G\in\G(w)\}$;
\item $\G(w)$ is admissible.
\end{enumerate}

To establish (i) and (ii), we apply the following lemma.

\begin{lemma}[see Schrijver \cite{Schrijver:1986}]\label{lem:duality}
Let $M$ be a rational $m$ by $n$ matrix and $c$ be a rational vector of length $m$.
Then $Mx=c$ has an integral solution if and only if for all rational vectors $y$ of length $n$,
$y^Tc$ is integral whenever $y^TM$ is integral.
\end{lemma}

For convenience, we write $a\equiv b$ if $a-b$ is an integer.

\begin{proof}[Proof of (i)]
For $i\in [q-1]$ 
let $X_i$ be rationals such that integral conditions hold.
In other words, for all $G\in\G(w)$, we have
\begin{equation*}
\sum_{i=1}^{q-1} X_im_i(G)\equiv 0.
\end{equation*}

For $i\in [q-1]$, consider the digraph 
$G(\w)$ with $w_i=w$ and $w_j=0$ for $j\ne i$.
Hence, we have
\begin{equation}\label{eq1}
w(w-1)X_i\equiv 0.
\end{equation}
For $(i,j)\in \overline{\binom{[q-1]}{2}}$, consider the digraph
$G(\w)$ with $w_i=w-1$, $w_j=1$ and $w_k=0$ for $k\notin\{i,j\}$.
\begin{equation}\label{eq2}
(w-1)^2X_i+(w-1)X_j\equiv 0.
\end{equation}
Subtracting \eqref{eq2} from \eqref{eq1},
we have $(w-1)X_i\equiv (w-1)X_j$ for all $i, j$.
Since $w-1$ divides $n-1$, we have
 
 \begin{equation*}
(n-1)X_i\equiv (n-1)X_j.
\end{equation*}

Finally, we compute
\begin{equation*}
n(n-1)\sum_{i=1}^{q-1}X_i \equiv (q-1)n(n-1)X_1\equiv 0,
\end{equation*}
since $w(w-1)$ divides $(q-1)n(n-1)$ and \eqref{eq1} holds.
Then (i) follows from Lemma \ref{lem:duality}.
\end{proof}

\begin{proof}[Proof of (ii)]
For $i\in [q-1]$, 
let $X_i$ and $Y_i$ be rationals such that integral conditions hold.
In other words, for all $G\in\G(w)$ and $u\in G$, we have
\begin{equation*}
\sum_{i=1}^{q-1} X_i\ind_i(u,G)+Y_i\outd_i(u,G)\equiv 0.
\end{equation*}

For $i\in [q-1]$, consider the digraph 
$G(\w)$ with $w_i=w$ and $w_j=0$ for $j\ne i$
and consider any vertex in $G(\w)$.
Hence, we have $(w-1)X_i+(w-1)Y_i\equiv 0$.
Since $(w-1)$ divides $(n-1)$,

\begin{equation*}
(n-1)X_i+(n-1)Y_i\equiv 0.
\end{equation*}

Therefore, the following relation is immediate.
\begin{equation*}
(n-1)\sum_{i=1}^{q-1}(X_i +Y_i) \equiv 0,
\end{equation*}

Again (ii) follows from Lemma \ref{lem:duality}.
\end{proof}

\begin{proof}[Proof of (iii)]
Summing up $\mu(G)$ for $G\in \G(w)$, 
we obtain a constant vector of length $q-1$ by symmetry.
Therefore, admissibility of $\G(w)$ is immediate.
\end{proof}

%
%
%
%
%
%
%

\bibliographystyle{IEEEtran}
\bibliography{mybibliography.bib}

\end{document}